\documentclass{amsart}

\usepackage{a4wide,latexsym,amssymb,amsthm,amsmath,color}

\theoremstyle{plain}

\newtheorem{theorem}{Theorem}

\newcommand{\E}{\mathbb{E}}
\newcommand{\C}{\mathbb{C}}

\newcommand{\me}{M^n\text{ in }\mathbb{E}^{n+m}}
\newcommand{\mc}{M^n\text{ in }\mathbb{C}^n}

\theoremstyle{remark}

\makeatletter
\def\Ddots{\mathinner{\mkern1mu\raise\p@
		\vbox{\kern7\p@\hbox{.}}\mkern2mu
		\raise4\p@\hbox{.}\mkern2mu\raise7\p@\hbox{.}\mkern1mu}}
\makeatother

\title[Extrinsic principal directions and curvatures]
{On the extrinsic principal directions and curvatures of Lagrangian submanifolds }

\author[Moruz]{Marilena Moruz}
\address{ Marilena Moruz\newline KU Leuven, Department of Mathematics, Section of Geometry, Celestijnenlaan 200B, 3001 Leuven, Belgium}
\email{marilena.moruz@kuleuven.be}
\author[Verstraelen]{Leopold Verstraelen}
\address{Leopold Verstraelen\newline PiT and CiT,\newline KU Leuven, Department of Mathematics, Section of Geometry, Celestijnenlaan 200B, 3001 Leuven, Belgium}
\email{leopold.verstraelen@kuleuven.be}

\begin{document}


\maketitle

\begin{abstract}
From the basic geometry of submanifolds will be recalled what are \emph{the extrinsic principal tangential directions}, (first studied by Camille Jordan in the $18$seventies), and what are  \emph{the principal first normal directions}, (first studied by Kostadin Tren\u cevski in the $19$nineties), and what are \emph{their corresponding Casorati curvatures}. For reasons of simplicity of exposition only, hereafter this will merely be done explicitly in the case of arbitrary submanifolds in Euclidean spaces. Then, \emph{for the special case of Lagrangian submanifolds} in complex Euclidean spaces, \emph{the natural relationships between these distinguished tangential and normal directions and their corresponding curvatures will be established.}
\end{abstract}

\section{The extrinsic tangential principal directions of submanifolds}

For general submanifolds $M^n$ of dimension $n \,(\geq 2)$ and of co-dimension $m\, (\geq 1)$ in Euclidean spaces $\mathbb{E}^{n+m}$, Jordan \cite{1} studied the extrinsic curvatures $c^{T}_u(p)$ at arbitrary points $p\in M$ in arbitrary tangential directions determined by vectors $u\in T_pM$, $\|u\|=1$. These are the curvatures $c^T_u(p)=(d\varphi_u/ds)^2(0)$, whereby $\varphi_u(s)\in [0,\Pi/2]$ denotes \emph{the angle in $\mathbb E^{n+m}$ between the tangent spaces $T_pM$ at $p$ and $T_qM$ at a nearby point $q\in M$ in the direction $u$ of $M$ at $p$}, $s$ being an arclength parameter of a curve $\gamma$ on $M$ from $p=\gamma(0)$ in the direction $u=\gamma'(0)$ to $q=\gamma(s)$. And, he defined \emph{the tangential principal curvatures } $c^T_1(p)\geq c^T_2(p)\geq\ldots\geq c^T_n(p)\geq 0 $ \emph{of a submanifold $M^n$ in $\mathbb{E}^{n+m}$ at $p$ as the critical values of the tangential Casorati curvature function at $p$}, that is of the function $c^T(p):S^{n-1}_p(1)=\{ u\in T_pM| \| u\|=1 \} \to \mathbb{R}^{+}: u\mapsto c^T_u(p)$, and, \emph{he defined the tangential principal directions of a submanifold $M^n$ in $\mathbb{E}^{n+m}$} at $p$ as the directions in which these critical values of the curvatures $c^T_u(p)$ are attained, and proved these directions to be mutually orthogonal, say to be determined by orthonormal vectors $f_1, f_2,\ldots, f_n\in T_pM$.\\
In the first step of his original fundamental studies of the geometry of submanifolds, Tren\u cevski \cite{2,3,4,5} re-considered this work of Jordan, and, later, Stefan Haesen and Daniel Kowalczyk and one of the authors \cite{6} basically re-did this once again. In the latter paper were followed \emph{the 1890 Casorati views on the intuitively most natural scalar valued curvatures ``as such" of surfaces $M^2$ in $\mathbb{E}^3$}; (and, in \cite{6} and in \cite{7}, some tangential and normal kinds of curvatures of Riemannian submanifolds were started to be named after Casorati). Accordingly, in \cite{6}, the above \emph{tangential Casorati curvatures} rather came up as $c^T_u(p)=(d\psi_u/ds)^2(0),$ whereby $\psi_u(s)$ denotes \emph{the angle in $\mathbb{E}^{n+m}$ between the normal spaces $T_p^{\perp}M$ at $p$ and $T_q^{\perp}M$ at a nearby point $q$ in the direction of $u$;} (as already was known by Jordan, $\psi_u=\varphi_u$). And, \emph{as shown by Tren\u cevski, the extrinsic principal unit tangential vector fields $F_1,F_2,\ldots,F_n$ of a submanifold $M^n$ in $\mathbb{E}^{n+m}$ and their corresponding tangential Casorati principal curvature functions $c_1^T, c_2^T, \ldots, c_n^T: M\to \mathbb{R}^+ : p\mapsto c_1^T(p),c_2^T(p), \ldots, c_n^T(p)$ essentially are the orthonormal eigen vectors fields and their corresponding eigen functions of the symmetric linear Casorati operator $A^C=\sum\limits_{\alpha}A_{\alpha}^2$, whereby $A_{\alpha}=A_{\xi_{\alpha}}$ are the  shape operators of $M^n$ in $\mathbb{E}^{n+m}$ for arbitrary orthonormal normal frame fields $\xi_1,\xi_2,\ldots,\xi_m$ on  $M^n$ in $\mathbb{E}^{n+m}$, such that $A^C F_i=c_i^T F_i,$ $i\in\{1,2,\ldots,n\}$, $ \alpha \in \{1,2\ldots,m\}$; } (the \emph{intrinsic} principal tangential directions and their corresponding curvatures of a submanifold  $M^n$ in $\mathbb{E}^{n+m}$, of course, being its \emph{Ricci} principal directions and curvatures).\\
From the above, in particular, one may notice that \emph{ for hypersurfaces  $M^n$ in $\mathbb{E}^{n+1}$ the extrinsic principal tangential directions are ``the classical" principal directions of these hypersurfaces, whereas $\{c_1^T, c_2^T,\ldots, c_n^T \}=\{k_1^2,k_2^2,\ldots,k_n^2\}$, $k_1,k_2,\ldots,k_n$ being the classical principal curvatures of these hypersurfaces} corresponding to Kronecker's extension of Euler's theory of the curvature of surfaces $M^2$ in $\mathbb E^3$ to hypersurfaces $M^n$ in $\mathbb E^{n+1}$ for all dimensions $n\geq 2$.

\section{Felice Casorati's study of surfaces $M^2$ in $\mathbb E^3$}

\emph{Casorati \cite{8} defined his extrinsic scalar valued curvature $C(p)$ of a surface $M^2$ in $\mathbb E^3$ at one of its points $p$ as follows.} On $M^2$, consider a small geodesic circle $\gamma_{\Delta_{\rho}}$ centered at $p$ with radius  $\Delta_{\rho}$. Let $q$ be any point on $\gamma_{\Delta_{\rho}}$ and consider the geodesic $\delta$ from $p$ to $q$
parametrised by arclength, such that $p=\delta(0)$ and $q=\delta(\Delta_{\rho})$; at $p$, this geodesic points in the tangential direction $\delta'(0)=u$ to $M^2$ at $p$. Let $\eta(p)$ and $\eta(q)$ be the unit normal vectors on the surfaces $M^2$ in $\mathbb E^{3}$ at $p$ and at $q$ respectively, corresponding to a choice of unit normal vector field $\eta$ around $p$ on $M^2$ in $\mathbb E^{3}$. Then, \emph{in Casorati's words, according to our common sense, the angle $\Delta\psi_u$ between $\eta(p)$ and $\eta(q)$ does measure well how much the surface $M^2$ at $p$ curves in the direction $u$; the more the surface  curves in the direction $u$, the larger this angle.} Then, joining all the points $\delta(\Delta\psi_u)$ that thus correspond to all the points $q$ on the geodesic circle $\gamma_{\Delta_{\rho}}$ around $p$, associated with $\gamma_{\Delta_{\rho}}$, one obtains  on $M^2$ a closed curve $\Gamma_{\Delta_{\rho}}$ (which actually passes through $p$ whenever at $p$ the surface is not curved at all in some tangential directions $u$).  And, hence, according to our common sense, the bigger or the smaller the area's $A(\Gamma_{\Delta_{\rho}})$ enclosed on $M^2$ by the curves $\Gamma_{\Delta_{\rho}}$  as compared to the area's $A(\gamma_{\Delta_{\rho}})$ of the geodesic discs on $M^2$ bounded by the geodesics $\gamma_{\Delta_{\rho}}$, the more or the less that the surface $M^2$ ``as such"  in $\mathbb{E}^3$ is curved at $p$. \emph{ It was along this line of thought that Casorati defined his curvature of a surface $M$ in $\mathbb{E}^3$ at $p$ as $C(p)=\underset{\Delta_{\rho}\to 0}{lim} (A(\Gamma_{\Delta_{\rho}})/A(\gamma_{\Delta_{\rho}}) ) $, and he proved that $C(p)=\frac{1}{2}tr A^2(p)=\frac{1}{2}(k_1^2+k_2^2)(p)$ $=\frac{1}{2}\|h\|^2(p)$, whereby $k_1$ and $k_2$ are Euler's principal curvatures, $A$ is the shape operator of $M^2$ corresponding to $\eta$ and $h$ is the second fundamental form of $M^2$ in $\mathbb E^{3}$}.\\
At this stage, it might be not amiss to add the following comment.\,\emph{In the definition of his curvature $C$, Casorati followed the common basic idea from the original geometrical definitions of the curvature $K$ of Gauss and of the mean curvature $H$ of Germain via ratio's of well chosen area's related to the surfaces $M^2$ in $\E^3$.} For $K(p)$, these ratio's concern regions on $M^2$ around $p$ and their corresponding spherical images, and, for $H(p)$, these ratio's are for discs centered at $p$ in $T_pM$ and for the portions of the corresponding circular cylinders perpendicular to $T_pM$ in between $T_pM$ and the surface $M^2$ in $\E^3$ itself. And, while for the curvatures of Germain and Gauss this lead to \emph{the first two elementary symmetric functions of $k_1$ and $k_2$,} $H = \frac{1}{2}tr A=\frac{1}{2}(k_1+k_2)$ and $K=det A=k_1k_2,$ Casorati's geometrical definition of his curvature yielding that $C=\frac{1}{2}tr A^2=\frac{1}{2}(k_1^2+k_2^2)$ lead to \emph{the third elementary symmetric function of Euler's principal curvatures}.

\section{The first normal principal directions of submanifolds}

\emph{Tren\u cevski determined the maximal possible dimensions of the osculating spaces of all orders for submanifolds $M^n$ in $\E^{n+m}$ and, moreover, in the related succesive normal spaces, also of all orders, determined appropriate orthonormal frames of principal normal vector fields and corresponding principal normal curvatures.} For our present purpose, it may suffice here to restrict within this grand theory to what is stated in Theorem 1 of \cite{7}: ``\emph{The first principal normal directions of a submanifold $M^n$ in $\E^{n+m}$ are the normal directions of $\me$ in which the normal Casorati curvatures of $M^n$ attain their $m_1$ (=dimension of the first normal space $N_1$) non-zero critical values."} \emph{The first normal space $N_1$ of $\me$ } is the subspace of the total normal space $T^{\perp}M$ of $\me$ given by $N_1=Im\, h=\{h(X,Y)| X,Y \in TM\}$, whereby \emph{$h$ is the second fundamental form of the submanifold $M$,} or, still, $N_1$ is the orthogonal complement in $T^{\perp}M$ of the subspace consisting of \emph{all normals} $\xi$ with vanishing shape operators $A_{\xi}$, or, equivalently, \emph{with vanishing normal Casorati curvature $c^{\perp}_{\xi}$; $N_1=\{\xi\in T^{\perp}M| A_{\xi}=0 \}^{\perp}$
$=\{ \xi\in T^{\perp}M| c^{\perp}_{\xi} =0 \}^{\perp}$}, such that \emph{the first osculating space of $\me$ is given by $TM\oplus N_1$ }.\\
The considerations of Casorati on surfaces $M^2$ in $\mathbb{E}^3$ that were recalled above can straightforwardly be taken over to general submanifolds $\me$, (and to general submanifolds $M^n$ in ambient general Riemannian spaces $\tilde M^{n+m}$, for that matter), cfr. \cite{6}. And, in \cite{6}, 
a.o. one may find the property that the Casorati curvature (as such) of a submanifold $\me$ equals the arithmetic mean of its tangential principal Casorati curvatures: $C= \frac{1}{n}\| h\|^2=\frac{1}{n}tr A^C=\frac{1}{n}\sum\limits_{\alpha}tr A^2_{\alpha}=\frac{1}{n}\sum\limits_{i}c_i^T$. Moreover, it seems not without interest to observe that $C_{\xi}(p)=\frac{1}{n}tr A^2_{\xi}(p)$ is \emph{the Casorati curvature (as such) at $p$ of the projection $M^n_{\xi p}$ of the submanifold $\me$ onto the $(n+1)D$ subspace $\E^{n+1}$ of $\E^{n+m}$ which is spanned by $T_pM=\mathbb R^n$ together with the normal line $[\xi(p)]$,} $\xi$ being any unit normal vector field on $\me$, and, hence, that  $C_{\xi}(p)=\frac{1}{n}\sum\limits_i c^{T}_{\xi i}(p)$, i.e. $C_{\xi}(p)$ is the arithmetic mean of the tangential Casorati curvatures $c^{T}_{\xi i}$ of this projected hypersurface $M^n_{\xi}$ at $p$; (for some general considerations relating the contemplation and the theory of submanifolds, see \cite{9}). The functions $c^{\perp}_{\xi}: \mathbb{S}^{m-1}(1)=\{\xi\in T^{\perp}M\,|\, \| \xi\|=1  \}\to \mathbb R^+: \xi\mapsto c^{\perp}_{\xi}=\frac{1}{n} tr A^2_{\xi}$, are called \emph{the normal Casorati curvatures of $\me$;} more precisely, \emph{the normal Casorati curvature of $\me$ in the direction determined by a unit normal vector field $\xi$ is defined as $c^{\perp}_{\xi}=\frac{1}{n}tr A^2_{\xi}$}.\\
In the total, $mD$ normal space $T^{\perp}M$ of $\me$, consider the following symmetric linear operator $a: T^{\perp}M\to T^{\perp} M: \xi\mapsto a(\zeta)=\frac{1}{n}\| \zeta\| \sum\limits_{\alpha} (tr A_{\zeta}A_{\alpha})\xi_{\alpha}; $ 
 (in \cite{10}, Bang-Yen Chen basically introduced this operator in the study of the submanifolds for which $a(\vec H)=\vec 0$, $\vec H$ being \emph{the mean curvature vector field} of $\me$, submanifolds which later were called \emph{Chen submanifolds}; in this respect, see also \cite{11} and \cite{12}). And, by the principal axes theorem, \emph{there exists an orthonormal frame} $\eta_1,\eta_2,\ldots,\eta_{m_1}, \eta_{m_1+1}, \ldots, \eta_{m}$ \emph{of eigen vector fields for this operator $a:T^{\perp}M\to T^{\perp}M $ ($m_1=dim N_1$), with corresponding eigen functions $c_1^{\perp}=\frac{1}{n}tr A^2_{\eta_1}\geq  c_2^{\perp}=\frac{1}{n}tr A^2_{\eta_2}\geq \ldots \geq c_{{m_1}}^{\perp}= \frac{1}{n}tr A^2_{\eta_{m_1}} $ $> c_{m_1+1}^{\perp}=tr A^2_{\eta_{m_1+1}}=\ldots=c_m^{\perp}=tr A^2_{\eta_m}=0.$} The normal vector fields $\eta_1,\eta_2,\ldots,\eta_{m_1}$ span the first normal space $N_1$ of $\me$ and, following Tren\u cevski, are called \emph{the first principal normal vector fields of the submanifold $\me$ with corresponding first principal normal curvatures $c_1^{\perp}\geq c_2^{\perp}\geq\ldots\geq c_{m_1}^{\perp}>0$.}  So, with indices $\alpha_1\in \{ 1,2,\ldots,m_1 \}$, $\{ \eta _{\alpha_1} \}$ \emph{is an orthonormal frame field of the first normal space $N_1$ for which $a(\eta_{\alpha_1}) =c_{\alpha_1}^{\perp}\eta_{\alpha_1}$, whereby $c_{\alpha_1}^{\perp}$ $=\frac{1}{n}tr A^2_{\alpha_1} (>0)$ are the principal normal Casorati curvatures of $\me$}.
 \noindent 
\section{The principal tangent and the first principal normal directions of Lagrangian submanifolds}
From \emph{Section 16: Totally real and Lagrangian submanifolds of K\"ahler manifolds} of Chen's contribution on \emph{Riemannian submanifolds} in \cite{11}, is taken the following: ``The study of totally real submanifolds of a K\"ahler manifold from differential geometric points of views was initiated in the early 1970's. ( -- By Bang-Yen Chen and Koichi Ogiue \cite{13} -- ; the authors.) A totally real submanifold $M$ of a K\"ahler manifold $\tilde M$ is a submanifold such that the almost complex structure $J$ of the ambient manifold $\tilde M$ carries each tangent space of $M$ into the corresponding normal space of $M$, that is, $J(T_pM)\subset T_p^{\perp}M$ for any point $p\in M$. ($\ldots$) \emph{ A totally real submanifold $M$ of a K\"ahler manifold $\tilde M$ is called Lagrangian if $dim_{\mathbb R}M=dim_{\mathbb C}\tilde M$}. $1$-dimensional submanifolds, that is, real curves, in a K\"ahler manifold are always totally real. For this reason, we only consider totally real submanifolds of dimension $\geq 2$.($\ldots$) \emph{For a Lagrangian submanifold $M$ of a K\"ahler manifold $(\tilde M, g, J)$ the tangent bundle $TM$ and the normal bundle $T^{\perp}M$ are isomorphic via the almost complex structure $J$ of the ambient manifold.} In particular, this implies that \emph{the Lagrangian submanifold has flat normal connection if and only if the submanifold is a flat Riemannian manifold.}"\\
To continue in our aim to go for simplicity and concreteness of presentation, (although, clearly, the following matters do hold more generally), next we do restrict our attention to \emph{the real $n$ dimensional totally real submanifolds $M^n$ of the complex $n$ dimensional complex Euclidean spaces} $\tilde M^n=\mathbb C^n=(\E^{2n}, \tilde J),$ that is, to \emph{the Lagrangian submanifolds $M^n$ in $\mathbb C^n$,} thus having $\tilde J (TM)=T^{\perp}M$ and $\tilde J (T^{\perp}M)=TM$, $\tilde J$ being \emph{the complex structure of the Kaehler manifold $\tilde M^n$}. On $M$ in $\tilde M$, \emph{tangential vector fields} will be denoted by $X,Y,Z,\ldots$ and \emph{normal vector fields} by $\xi, \eta, \zeta,\ldots$ . Further, let $\tilde g$ and $\tilde \nabla$, respectively $g$ and $\nabla$, be \emph{the metrics and the corresponding Riemannian connections} on $\tilde M$ and $M$, respectively.  The \emph{equations of Gauss and of Weingarten} are given by
\begin{align}
\tilde\nabla_XY=\nabla_XY+h(X,Y),\label{1}\\
\tilde{\nabla}_X\xi=-A_{\xi}(X)+\nabla^{\perp}_X\xi,\label{2}
\end{align}
whereby $\nabla^{\perp}$ is \emph{the normal connection} and $h$ \emph{the second fundamental form} and \emph{$A_{\xi}$ the shape operator with respect to} $\xi$ of the submanifold $M$ in $\tilde M$, so that
\begin{equation}\label{3}
\tilde g(h(X,Y),\xi)=g(A_{\xi}(X), Y).
\end{equation}
Applying the complex structure $\tilde J$ to \eqref{1}, it follows that 
\begin{equation}\label{4}
\tilde J(\tilde \nabla_XY)=\tilde J(\nabla_XY)+\tilde J(h(X,Y)),
\end{equation}
while writing \eqref{2} out for $\xi=\tilde J Y$, it follows that 
\begin{equation}\label{5}
\tilde{\nabla}_X(\tilde JY)=-A_{\tilde JY}(X)+\nabla_X^{\perp}(\tilde JY).
\end{equation}
By the parallelity of $\tilde J$, $\tilde \nabla\tilde J=0$, or, still, $\tilde J(\tilde \nabla_XY)=\tilde{\nabla}_X(\tilde JY)$, the left-hand sides in \eqref{4} and\eqref{5} are equal, and, hence, in particular, also the tangential components of the right-hand sides in \eqref{4} and \eqref{5} are equal:
\begin{equation}\label{6}
\tilde J(h(X,Y))=-A_{\tilde JY}(X).
\end{equation}
Writing out \eqref{3} for $\xi=\tilde JZ$, it follows that 
\begin{equation}\label{7}
\tilde g(h(X,Y), \tilde JZ)=g(A_{\tilde JZ}(X), Y),
\end{equation}
which by \eqref{6} leads to
\begin{equation}\label{8}
\tilde g( h(X,Y), \tilde JZ)=g(-\tilde J(h(X,Z)),Y).
\end{equation}
Since $\tilde J$ is almost complex, $\tilde J^2=-I$, and since \emph{$\tilde g$ is Hermitian,} so that $\tilde g(\tilde J\tilde V, \tilde J\tilde W)=\tilde g(\tilde V,\tilde W)$ for all vector fields $\tilde V$ and $\tilde W$, and, hence, in particular, for $\tilde V=-\tilde J(h(X,Z))$ and $\tilde W=Y$, \eqref{8} becomes
\begin{align}\label{9}
\begin{split}
\tilde g( h(X,Y), \tilde JZ)=&g(-\tilde J^2(h(X,Z)), \tilde JY )\\
=&g(h(X,Z), \tilde JY).
\end{split}
\end{align}
In view of its crucial importance in what comes next, we have cared to work out in detail this property from \cite{13} as obtained in \eqref{9}, which may be stated as follows.\emph{ For all tangential vector fields $X,Y,Z$ on a Lagrangian submanifold $M^n\subset\mathbb{C}^n \, (\tilde M^n)$}:
\begin{equation}\label{10}
\tilde g(h(X,Y), \tilde JZ)=\tilde g(h(X,Z),\tilde JY)=\tilde g(h(Y,Z),\tilde JX).
\end{equation}
 \emph{For any tangential orthonormal frame field $\mathcal{F}=\{E_1, E_2,\ldots, E_n  \}$ on a Lagrangian submanifold $M^n$}, $\tilde{ \mathcal{F}}=\{E_1, E_2,\ldots, E_n, \xi_1=\tilde{J}E_1,\xi_2=\tilde{J}E_2, \ldots, \xi_n=\tilde{J}E_n  \}=  \{ E_i, \xi_i=\tilde{J} E_i\} $, $(i,j,k,\alpha, \beta\in \{1,2,\ldots,n\}  )$ \emph{is a corresponding adapted orthonormal frame field of $\mathbb{C}^n (\tilde M^n)$ along $M^n$.  The local coordinates of the operator $A^C:TM\to TM $ of Casorati and of the operator $a:T^{\perp}M\to T^{\perp}M$ of Tren\u cevski with respect to such frame fields $\tilde{\mathcal{F}}$ are given by } 
\begin{align}\label{11}
\begin{split}
A^C_{ik}=&(\sum\limits_{\alpha}A_{\alpha}^2 )_{ik}=\sum\limits_{\alpha}(A_{\alpha}^2)_{ik} \\
=&\sum\limits_{\alpha}\sum\limits_{j} h^{\alpha}_{ij}  h^{\alpha}_{jk}
\end{split}
\end{align}
and 
\begin{align}\label{12}
\begin{split}
a_{\alpha \beta}=&tr(A_{\alpha}A_{\beta})\\
=&\sum\limits_i\sum\limits_j h_{ij}^{\alpha} h_{ji}^{\beta},
\end{split}
\end{align}
whereby $h_{ij}^{\beta}$ are the local coordinates of the symmetric second fundamental form $h: TM\times TM\to T^{\perp}M$. So, $E_i$ \emph{determines a Casorati principal tangential vector field on $M^n$ in $\mathbb{C}^n$ with corresponding principal tangential Casorati curvature $c_i^T$ if and only if }
\begin{equation}\label{13}
\forall k\neq i\ : \ A^C_{ik}=\sum\limits_{\alpha}\sum\limits_{j}h_{ij}^{\alpha} h^{\alpha}_{jk}=0,
\end{equation} 
whereby then
\begin{equation}\label{14}
c_i^T=A^C_{ii}=\sum\limits_{\alpha}\sum\limits_{j}(h_{ij}^{\alpha})^2,
\end{equation}
and, $\xi_{\alpha}=\tilde J E_{\alpha}$ determines a first principal normal vector field, or, first Casorati principal normal vector field (as these vector fields later on also might be termed), \emph{on $M^n$ in $\mathbb{C}^n$ with corresponding principal normal Casorati curvature $c^T_{\alpha}$ if and only if } 
\begin{equation}\label{15}
\forall \beta\neq\alpha\ : \ a_{\alpha\beta}=\sum\limits_i\sum\limits_j h_{ij}^{\alpha} h_{ji}^{\beta}=0,
\end{equation}
\emph{whereby then}
\begin{equation}\label{16}
c_{\alpha}^{\perp}=a_{\alpha\alpha}=\sum\limits_i\sum\limits_j (h_{ij}^{\alpha})^2.
\end{equation}
Written in local coordinates, \emph{the above property \eqref{10} amounts to}
\begin{equation}\label{17}
\forall i,j,k\ :\ h_{ij}^k=h_{ik}^j=h_{jk}^i,
\end{equation}
so that, from \eqref{13} and \eqref{15}, and, from \eqref{14} and \eqref{16}, in particular, we may conclude the following.
\begin{theorem}\label{T1}
Let $M^n$ be a Lagrangian submanifold of the complex Euclidean space $\mathbb{C}^n$ (or, of any Kaehler manifold $\tilde M^n$). Then, a tangential vector field $T$ is a tangential principal Casorati vector field with corresponding tangential Casorati principal curvature $c^T(>0)$ if and only if $N=\tilde J T$ is a normal principal Casorati vector field -- whereby $\tilde J$ is the complex structure of $\mathbb C^n$ (or, of the ambient Kaehler space  $\tilde M^n$) -- with corresponding normal Casorati principal curvature $c^{\perp}=c^T(>0)$.
\end{theorem}
	
\begin{theorem}\label{T2}
		Let $M^n$ be a Lagrangian submanifold of the complex Euclidean  space $\mathbb C^n$ (or, of any ambient Kaehler manifold $\tilde M^n$) with first normal space of maximal dimension ($m_1=dim N_1=n=co-dim M$). Then, $M^n$ admits an adapted orthonormal frame field $\tilde{ \mathcal{F}}=\{ F_1, F_2, \ldots, F_n, \eta_1=\tilde J F_1, \eta_2=\tilde J F_2,\ldots, \eta_n=\tilde J F_n  \} $ in $\C^n\ (\tilde M^n)$ of which the $n$ tangential vector fields are the principal Casorati tangential vector fields and of which the $n$ normal vector fields are the principal Casorati normal vector fields of $\mc$ $(\tilde M^n)$, and the corresponding tangential and normal principal curvatures are equal, $(\forall i\ :\ c_i^T=c_i^{\perp})$.
\end{theorem}


\noindent {\bf Acknowledgements.}
The author named first is a Postdoctoral Fellow of \emph{The Research Foundation -- Flanders (FWO)}.

\end{document}